\documentclass[12pt, reqno]{amsart}
\usepackage{amsmath, amstext, amsbsy, amssymb, amscd}
\usepackage{amsmath}
\usepackage{amsxtra}
\usepackage{amscd}
\usepackage{amsthm}
\usepackage{amsfonts}
\usepackage{amssymb}
\usepackage{eucal}
\usepackage{color}
\usepackage{multicol}
\usepackage{multirow}

\newcommand{\ds}{\displaystyle}

\begin{document}

\title[Zombie Dice: An Optimal Play Strategy]
{Zombie Dice: An Optimal Play Strategy}

\author{Heather L. Cook}
\address{Department of Mathematics, Computer Science, and Physics, Roanoke College, Salem, VA 24153}
\email{hlcook@mail.roanoke.edu}

\author{David G. Taylor}
\address{Department of Mathematics, Computer Science, and Physics, Roanoke College, Salem, VA 24153}
\email{taylor@roanoke.edu}

\begin{abstract}
We discuss the game of Zombie Dice, published by Steve Jackson Games. This game includes green, yellow, and red dice. Each die has brain, footprint, and shotgun symbols on it, with each color of dice having a different amount of each symbol. Out of the dice, three are randomly picked and rolled. The player plays as if he were the zombie, meaning that brains are wanted and shotguns are not. Footprints are rerolled if the player chooses to keep going and not score. One brain equals one point. If three shotguns are accumulated, then that player's turn is over and he loses all his brains for no points (busting). The objective of the game is to gain thirteen or more points. In this article, we investigate a model for deciding whether or not to continue rolling (if given the opportunity). With this model, we create a decision point given information about the current player's turn including the amount and color of dice left in the cup, color and number of footprints, and the current number of shotguns of the player. Examples will be shown to highlight the game and its strategy fashioned from the model.
\end{abstract}

\maketitle
\section{History}

People from all over the world have been playing games for centuries with an early form of dice being frequently discovered in Assyrian and Simerian archeological digs. A man named Girolamo Cardano $(1501-1576)$ is sometimes considered the first to give a working definition for probability with his works on this not being published until close to a century after his passing, around $1676$. By this time though, others had begun writing about the mathematics behind games including Galileo specifically around $1613$ and $1623$, who especially focused on those involving dice. Galileo had noted in \textit{Concerning an Investigation of Dice} the amount of sums when given two or more dice and how many ways this could occur in addition. Antoine Gombaud, also known as The Chevalier de M\'{e}r\'{e}, $(1607-1684)$ became caught up in gambling in which he began to question explanations to certain aspects that are naturally mathematical. Specifically, he questioned rolling at least one six on a die in four rolls with an even money bet in which he had gained a considerable sum of money while, over time, he had noticed that he lost more money than he expected while betting on one or more double sixes in twenty-four rolls of two dice and wondered why. In $1654$, he asked for help from Blaise Pascal, the well-known mathematician, in order to understand. Sparking Pascal's interest, Pascal corresponded with Pierre de Fermat to answer de M\'{e}r\'{e}'s question in which probability was formally stated.

Using the Multiplication Principle, which is the product of the number of ways each task can be performed, we can explain de M\'{e}r\'{e}'s dice problem. To start with the game of obtaining at least one six in four rolls, there are $6^4=1296$ possible outcomes. Now, we need to account for how many of these actually contain at least one six. This is done by counting how many of these outcomes do \emph{not} contain a six. So, if one of the dice has to be a six, then we have five outcomes per one roll which gives $5^4=625$ outcomes that do not have a six. So, taking the first minus the second number of outcomes we have $1296-625=671$ with at least one six. Then the proportion of outcomes for at least one six with the four rolls is
$$ \frac{671}{1296}=.5178.$$ We can see from this that we should expect to win more often than not.
With the second game of one or more double sixes in 24 rolls with two dice, however, the sample space is comprised of $6^2=36$ separate outcomes. Within this sample space, there is only one pair of double sixes which means $35$ of these outcomes we don't want. So, the proportion of these outcomes that do not contain double sixes is
 $$ \left(\frac{35}{36}\right)^{24}=.5086. $$
Since this is greater than one half, we know that we should expect to lose more often. Therefore, we can see the solution to de M\'{e}r\'{e}'s queries. This is all how probability started, leading to everything we know today about it.

\section{Introduction to Zombie Dice}

The game Zombie Dice was created by Steve Jackson Games, Inc. In the game of Zombie Dice, the player plays as if he is the zombie, eating brains and getting shotgunned. For each die, there are three different face types. First are the brains which means eating the victim's brain, and they are set aside for that player's round. Shotguns are the zombie, or player, being ``shot''. These dice also are set aside. The last face type are footprints. This is where the victim escaped, so these dice are kept in front to roll again if the player so chooses.

There are thirteen dice total with three red, four yellow, and six green. Each specific color of dice have separate amounts of brains and shotguns, but each have two footprints. The red dice have one brain and three shotguns, the yellow have two brains and two shotguns, and the green have three brains and one shotgun.

The player starts by shaking the cup of thirteen dice and picking three random dice. If the player rolls shotguns or brains, they are placed to the side until the end of their turn. The player's turn is over when at least three shotguns are rolled in total for that round in which he receives no points for that round. The player can also choose to stop and score his points instead of rolling again, in which he gains one point per one brain accumulated. When a player chooses to continue, he then takes enough random dice from the cup to total three, which includes any footprints from the previous roll. For instance, say the player had one footprint from the previous roll, then only two dice are taken from the cup for a total of three dice to roll. If there are not any dice left in the cup, or not enough to total three for the next roll, the player notes how many brains he has gained and keeps the shotguns out while returning the rest of the dice to the cup and continues as normal. The object of the game is to gain thirteen or more points total. Once a player has obtained this, the rest of the players in that round finish. Whoever has the most brains at the end of the round wins. A tiebreaker round is played if a tie occurs.

The question we analyzed was, given we currently have $X$ brains, should we keep rolling or score our points? We began by creating two functions in order to find the probability of getting three or more shotguns given our current number of shotguns and the probability of getting zero, one, two, or three brains on the next roll which gave the following probability distribution functions
\begin{multicols}{2}
\begin{center}\begin{tabular}{c|c} $Y$ & $P(Y)$ \\\hline 0 & 0.347449 \\ 1 & 0.444833 \\ 2 & 0.183372 \\ 3 & 0.024346 \end{tabular} \end{center}
\columnbreak
\begin{center}\begin{tabular}{c|c} $X$ & $P(X)$ \\\hline 0 & 0.245144 \\ 1 & 0.444056 \\ 2 & 0.261072 \\ 3 & 0.049728 \end{tabular} \end{center}
\end{multicols} \noindent where $Y$ is the number of shotguns and $X$ is the number of brains. We then created a model which led to a recursive program in order to obtain the expected values for all future rolls. With this, we created a strategy along with generalizations to a simpler strategy that's easily memorizable. In order to reach the goal of optimal strategy to help answer the question though, we need a reminder and the definition for expected value. \emph{Expected value} can be used to see how much one should expect to receive over repeated trials, or in the long run. To give a more formal definition, we need to know about random variables. With an experiment having a sample space $S$, a function $X$ that assigns only a single real number to each element of $S$, $X(s)=x$, is called a \emph{random variable}. Now to define expected value, we let $X$ be a random variable and $e_1, e_2, \ldots, e_k$ be the elements of the sample space with the corresponding probabilities $p_1, p_2, \ldots, p_k$. Then the expected value $(EV)$ of the random variable is
$$EV(X)=p_1\cdot X(e_1)+p_2\cdot X(e_2)+\ldots+p_k\cdot X(e_k).$$ Now we can look at conditional expectation. First, we need to let $X$ be a discrete random variable and $E$ to be an event where $P(E)>0$. Then the \emph{conditional expectation} of the random variable $X$ given $E$ has occurred is defined as
$$EV(X|E)=\sum_{x\in X}xP(X=x|E).$$

Our answer to the question posed above is that the player will continue to roll if
$$- X \cdot\ P(B|S=Y) + EV \cdot\ (1-P(B|S=Y)) > 0$$
where
$$\begin{aligned} X &= \text{current number of brains} \\ Y &= \text{current number of shotguns} \\ P(B|S=Y) &= \text{probability of busting given we have $Y$ shotguns} \\ EV &= \text{expected number of brains on next rolls.}\end{aligned}$$

\section{The Model}

In order to construct one recursive function to answer when the player should keep rolling or stop and score points, we use the multinomial distribution along with summations. The multinomial coefficient counts the number of ways to choose from $n$ objects a collection of $k$ types of objects where there are $x_1$ of one type, $x_2$ of another, and so on up to $x_k$ of the final kind. We require that $x_1+x_2+\ldots+x_k=n$ so each $n$ object will appear as one of the types. The multinomial coefficient is given by
$$\ds{n\choose x_{1},x_{2},...,x_{k}}=\frac{n!}{x_1!\cdot x_2!\cdots x_k!}.$$
The multinomial distribution uses the multinomial coefficient and has $n$ trials which can result in any number $k$ of different results. In our case, the different results are brains, footprints, and shotguns. The trials are identical so the probability of outcome $R_{i}$ stays the same for each repetition. If the probability of outcome $R_{i}$ is given as $p_{i}$ then the probability that the experiment yields an $x_{i}$ number of outcome $R_{i}$ is given by:
$$\text{P}(R_{1}=x_{1},R_{2}=x_{2},...,R_{k}=x_{k})=\ds{n\choose x_{1},x_{2},...,x_{k}}\cdot p_{1}^{x_{1}}\cdot p_{2}^{x_{2}}\cdot\cdot\cdot p_{k}^{x_{k}}$$
where we require that
$$\begin{aligned} x_{1}+x_{2}+...+x_{k}&=n \text{ and} \\
                  p_{1}+p_{2}+...+p_{k}&=1.\end{aligned}$$
With this definition we can now write our probabilities much easier. For example, we can define
$$\text{P(Green Shotguns)}=P_{\text{green}}^{(0,0,g)}=\underbrace{\ds{g\choose 0,0,g}}_{\text{all green shotguns}}\cdot \underbrace{\left(\frac{1}{2}\right)^0}_{\text{brains}}\cdot \underbrace{\left(\frac{1}{3}\right)^0}_{\text{footprints}}\cdot \underbrace{\left(\frac{1}{6}\right)^g}_{\text{shotguns}}.$$
This will now allow us to write the probability of getting three shotguns of some combination between the green, yellow, or red dice as
$$\text{P(3 Shotguns With }\{g,y,r\})=P_{\text{green}}^{(0,0,g)}\cdot P_{\text{yellow}}^{(0,0,y)}\cdot P_{\text{red}}^{(0,0,r)}.$$
The probability of obtaining a certain combination of dice from the cup, the coefficient on these probabilities, is is computed as
$$\text{P}(g\text{ Green, }y\text{ Yellow, }r\text{ Red Dice})=\frac{\overbrace{\ds{6\choose g}}^{\text{green dice}}\cdot \overbrace{\ds{4\choose y}}^{\text{yellow dice}}\cdot \overbrace{\ds{3\choose r}}^{\text{red dice}}}{\underbrace{\ds{13\choose 3}}_{\text{all dice combinations}}}.$$
In general, we can define the following which is the probability of getting $x_{b}$ brains, $x_{f}$ footprints, and $x_{s}$ shotguns for the green dice where we require that $x_{b}+x_{f}+x_{s}=n$ as
$$P_{\text{green}}^{(x_{b},x_{f},x_{s})}=\underbrace{\ds{n\choose x_{b},x_{f},x_{s}}}_{\text{orderings}}\cdot \underbrace{\left(\frac{1}{2}\right)^{x_{b}}}_{\text{brains}}\cdot \underbrace{\left(\frac{1}{3}\right)^{x_{f}}}_{\text{footprints}}\cdot \underbrace{\left(\frac{1}{6}\right)^{x_{s}}}_{\text{shotguns}}.$$
The formulas for the yellow and red dice are similar with the same restriction and can be written respectively as
$$P_{\text{yellow}}^{(x_{b},x_{f},x_{s})}=\ds{n\choose x_{b},x_{f},x_{s}}\cdot \left(\frac{1}{3}\right)^{x_{b}}\cdot \left(\frac{1}{3}\right)^{x_{f}}\cdot \left(\frac{1}{3}\right)^{x_{s}}$$
and
$$P_{\text{red}}^{(x_{b},x_{f},x_{s})}=\ds{n\choose x_{b},x_{f},x_{s}}\cdot \left(\frac{1}{6}\right)^{x_{b}}\cdot \left(\frac{1}{3}\right)^{x_{f}}\cdot \left(\frac{1}{2}\right)^{x_{s}}.$$
Note that the formula for the yellow dice can be simplified as
$$P_{\text{yellow}}^{(x_{b},x_{f},x_{s})}=\ds{n\choose x_{b},x_{f},x_{s}}\cdot \left(\frac{1}{3}\right)^{n}$$
since the probabilities for brains, footprints, and shotguns are the same.
The expected value of continuing is computed as
$$\text{Expected Value of Continuing}=\underbrace{-b\cdot PE}_{\text{round ends}}+\underbrace{EB\cdot (1-PE)}_{\text{round continues}}>0$$
meaning that, with this inequality, we can play the game as though we are betting with the current amount of brains being the wager and continue until the expected value becomes zero or negative.

As an example, using all of the above pieces together and defining $S$ as the set of all triples of the form $\{g,y,r\}$ where $0 \leq g,y,r \leq 3$ and $g+y+r=3$, we can compute the probability of gaining $3$ shotguns on the first round with the use of summations as
$$ \sum_{\{g,y,r\}\in S} \frac{\ds{6\choose g}\cdot \ds{4\choose y}\cdot \ds{3\choose r}}{\ds{13\choose 3}}\cdot p_{\text{green}}^{(0,0,g)} \cdot p_{\text{yellow}}^{(0,0,y)} \cdot p_{\text{red}}^{(0,0,r)}=\frac{94}{3,861}\approx0.02435$$
meaning that the chance of your turn ending is about $2.4\%$ on the first roll.

In order to generalize this formula even more, we can define the following $\mathbf{C}\{g,y,r\}$ where $g_c$ green, $y_c$ yellow, and $r_c$ red dice are left in the cup and $g_f$ green, $y_f$ yellow, and $r_f$ red footprints roll from the previous roll (the boldface reminds us that the functions depends on $g_c$, $y_c$, $r_c$, $g_f$, $y_f$, and $r_f$). It gives the probability of getting the combination $\{g,y,r\}$ from the cup. So,

$$C_{g_f,y_f,r_f}^{g_c,y_c,r_c}\{g,y,r\}=\mathbf{C}\{g,y,r\}=\frac{\ds{g_c\choose g-g_f}\cdot \ds{y_c\choose y-y_f}\cdot \ds{r_c\choose r-r_f}}{\ds{g_c+y_c+r_c\choose 3-(g_f+y_f+r_f)}}$$
where $\mathbf{C}\{g,y,r\}=0$ for impossible situations such as when $g_f=2$ with the combination $\{1,1,1\}$ since this would mean the player would reroll two green dice and therefore could not have a red \textit{and} a yellow die to roll as well.
This allows us to write the following which gives the probability of getting $3$ shotguns on the next roll as
$$P(\text{3 Shotguns})=\sum_{\{g,y,r\}\in S} \mathbf{C}\{g,y,r\}\cdot p_{\text{green}}^{(0,0,g)} \cdot p_{\text{yellow}}^{(0,0,y)} \cdot p_{\text{red}}^{(0,0,r)}.$$

We can tackle the expected amount of brains in a similar manner using the same notation as above. The expected number of brains on the next roll is then
$$\text{Expected Brains}=0\cdot \mathbf{B}(0)+1 \cdot \mathbf{B}(1)+2\cdot \mathbf{B}(2)+3\cdot \mathbf{B}(3)$$
where $\mathbf{B}(x)$ is the probability of obtaining $x$ brains on the next roll (which depends on the number and color of footprints and dice left in the cup).
Using the summation from the shotgun example above, we gain the following brain summation given the current values of $g_c$ and $g_f$ and the corresponding red and yellow dice in the cup and footprints. We write
$$\mathbf{B}(3)= \sum_{\{g,y,r\}\in S} \mathbf{C}\{g,y,r\}\cdot p_{\text{green}}^{(g,0,0)} \cdot p_{\text{yellow}}^{(y,0,0)} \cdot p_{\text{red}}^{(r,0,0)}.$$
There is a difference between the notation $\{g,y,r\}$ and $(i,j,k)$ though, in which the first represents a specific selection of dice where the player has $g$ green dice, $y$ yellow dice, and $r$ red dice to roll and the latter represents a specific color of dice getting $i$ brains, $j$ footprints, and $k$ shotguns. Each element in $S$, which is the set of all triples of the form $\{g,y,r\}$ where $0 \leq g,y,r \leq 3$ and $g+y+r=3$, gives a certain set of $(i,j,k)$ possibilities for each color separately. To make this easier, we denote $T(\{g,y,r\};x)$ to be the set of possibilities given values for $g, y, r, \text{and } x$. So the most general form for $\mathbf{B}(x)$ for other values of $x$ is then
$$\mathbf{B}(x)= \sum_{\{g,y,r\}\in S} \mathbf{C}\{g,y,r\}\sum_{T(\{g,y,r\};x)} p_{\text{green}}^{(g_i,g_j,g_k)} \cdot p_{\text{yellow}}^{(y_i,y_j,y_k)} \cdot p_{\text{red}}^{(r_i,r_j,r_k)}.$$
Then, we can construct a formula with summation notation to find the expected number of brains on the next roll given the dice that is left in the cup and the footprints to reroll as
$$\text{Expected Brains}=EB_{g_f,y_f,r_f}^{g_c,y_c,r_c}=\mathbf{EB}=\sum_{i=0}^3 i\cdot \mathbf{B}(i).$$

With this notation we can also denote the probability of getting $x$ shotguns given the regular parameters by $S_{g_f,y_f,r_f}^{g_c,y_c,r_c}=\mathbf{S}(x)$ which gives
$$\mathbf{S}(x)=\sum_{\{g,y,r\}\in S}\mathbf{C}\{g,y,r\}\sum_{T(\{g,y,r\};x)}p_{\text{green}}^{(g_k,g_j,g_i)} \cdot p_{\text{yellow}}^{(y_k,y_j,y_i)} \cdot p_{\text{red}}^{(r_k,r_j,r_i)}.$$
This is almost the same as that of the brains but with the difference of the use of $(g_k,g_j,g_i)$ where the $g_i$ and $g_k$ roles are reversed in $p_{\text{green}}^{(g_k,g_j,g_i)}$ since the order of $(g_k,g_j,g_i)$ corresponds to $g_k$ brains, $g_j$ footprints, and $g_i$ shotguns which is the reverse of the elements of $T$ when we deal with shotguns. The same goes for the corresponding red and yellow dice functions.

The probability that the round ends depends upon the number of current shotguns, so we need to address the probability of the round ending on the next roll due to shotguns. Since the player can only roll three shotguns maximum, the probability is as follows
$$\text{P(Round Ends)}=PE^{g_c,y_c,r_c}_{g_f,y_f,r_f}(s)=\mathbf{PE}(s)=\sum^{3}_{i=3-s}\mathbf{S}(i)$$
where we have $s$ current shotguns and the usual variables. If we have one current shotgun, then we need to add together the probabilities of attaining two or three shotguns on the next roll which is written as
$$\mathbf{PE}(1)=\mathbf{S}(2)+\mathbf{S}(3).$$
Now, we can rewrite the expected value of continuing with this new notation as
$$\text{Expected Value of Continuing}=\underbrace{-b\cdot\mathbf{PE}(s)}_{\text{round ends}}+\underbrace{\mathbf{EB}\cdot(1-\mathbf{PE}(s))}_{\text{round continues}}>0.$$
However, we need to correct this since the function $\mathbf{EB}$ gives brains that are not allowed. For example, if we have \{1,1,1\} as our element of $S$ that we are iterating over and we have two current shotguns, then the $T(\{1,1,1\};2)$ element has $(1,0,0)$ for the green dice, $(1,0,0)$ for the yellow dice, and $(0,0,1)$ for the red dice giving us two extra brains and ending the round since we also gained a third shotgun. Thus, some of the elements of $T$ need to be ignored for certain $S$ elements. Also, if the player chooses to reroll and does not obtain three shotguns, then he has the opportunity to roll again and get more brains. So, we need to consider the expected number of brains on all future rolls and not just the next roll. This is done recursively, easily through using computers.

Solving the Expected Value of Continuing equation for $b$, we can give a value of $b$ that is similar to a decision point for a fixed value of $s$ and the values for the quantity of dice remaining in the cup and current footprints. Thus, we get
$$\frac{\mathbf{EB}\cdot(1-\mathbf{PE}(s))}{\mathbf{PE}(s)}>b.$$
This means that the player should continue rolling as long as the number of brains the player has currently, $b$, is less than the value on the left-hand side after one roll. Each time the player rolls, the left-hand side needs to be recomputed since the expected value of brains and the probability of the round ending changes with each roll.

A program was written by Dr. David Taylor in which the limitations above were removed in order to obtain an optimal strategy for playing Zombie Dice. This program gives the decision point $b$ for each possible combination of the variables for the amount of dice in the cup, footprints, and current number of shotguns. For example, with the following set up
\begin{center}
\begin{tabular}{c|c|c|c|c}
Dice in Cup & Current FPs & Decision & Decision & Decision
\\ R Y G & R Y G & (SG 0) & (SG 1) & (SG 2)
\\\hline 2 3 1 & 1 0 1 & 78.338580 & 4.043669 & 0.180008
\end{tabular}
\end{center}
the player should continue rolling if he has less than or equal to $78$ brains with no current shotguns, or if he has less than or equal to $4$ brains with one current shotgun. If the player has two current shotguns, then he should stop rolling if he has one or more brains and he should only keep rolling if he has no brains because the left-hand side would have to be less than $0.180008$ to keep rolling. This would be the strategy when the player has $2$ red, $3$ yellow, and $1$ green dice in the cup as well as $1$ red, $0$ yellow, and $1$ green footprints to reroll. With every combination of colors and amounts of footprints and dice in the cup, the decision point changes and especially depends on the number of current shotguns.

\section{Strategy}

Using the table and strategy created through recursion as explained above, the player can be a formidable opponent throughout the entire game, especially during the beginning and middle. Near the end of the game however (which is when someone in the current round who is before the player in turn order reaches $13$ or more brains, or when an opponent who comes after the player gets $10$ or more brains) the player will want to reach $13$ or more brains no matter the cost in order to have a chance at winning. Since the expected number of brains with $3$ red, $4$ yellow, and $6$ green dice in the cup, no footprints, and no current shotguns (the first roll of a player's turn), is $3.315559$, then we know that if a player has $10$ dice from previous rounds and has yet to go, then you can expect them to attain at least $3$ more brains on their first roll for a total of $13$ to end the round without you having the chance to beat them and win. This is why a player would want to try and obtain as many points if someone before them has already gained $10$ points or more, which is deemed as a strategy for the end of the game.

There are many generalizations that can be made with this immense data set that is given by the written computer code which gives all possible combinations ($4867$ total) of the six usual parameters (the amount of red, yellow, and green dice left in the cup; and the amount of red, yellow, and green footprints) plus the number of current shotguns. Below is the table of the simple strategy sorted by the current number of shotguns which is slightly less optimal than using the enormous data set, but it is vastly easier.
  \begin{center}
        \begin{tabular}{|c|c|} \hline Current Shotguns & Rule \\ \hline\hline 0 & Always keep rolling. \\\hline
             \multirow{5}{*}{1} & If we must roll 3 red dice, stop at 1 brain. \\
              & If ($r_f = 2$ and $y_f=1$) or $y_c>g_c$, stop at 2 brains. \\
              & If ($r_f = 2$ and $g_f=1$) or $g_c>y_c$, stop at 3 brains. \\
              & If we must roll 3 green dice, roll again. \\
              & If $g_f = 2$, roll again. \\\hline
             \multirow{2}{*}{2} & With $g_f = 3$, stop at 2 brains. \\
              & Otherwise, stop at 1 brain. \\\hline
        \end{tabular}
    \end{center}

\section{Further Questions}

Throughout this paper, we have thoroughly explored the game of Zombie Dice and how one should optimally play. There may be other ways to present general strategies besides the way described above. We are certain that our strategy, especially the results from the recursive program, will beat the AI code for the iPad application when compared to see which plays more optimally.  We would like to thank Steve Jackson Games for sharing the iPad decision-making strategy with us (under an NDA) so that we will be able to do this comparison.

\end{document}